\documentclass[12pt]{amsart}
\usepackage{amsfonts,amssymb,amscd,amsmath,amsthm,enumerate,verbatim,calc}

\usepackage{times} 
\usepackage{mathptm}

\theoremstyle{plain}
\newtheorem{Theorem}{Theorem}[section]

\newtheorem{Lemma}[Theorem]{Lemma}
\newtheorem{Corollary}[Theorem]{Corollary}

\theoremstyle{definition}
\newtheorem{Definition}[Theorem]{Definition}

\newtheorem{Remark}[Theorem]{Remark}

\theoremstyle{remark}

\newtheorem*{chunk*}{}

\numberwithin{equation}{Theorem}

\newcommand{\mlabel}[1]%
  {\mbox{}\marginpar{\raggedleft\hspace{0pt}{\rm\ttfamily#1}}\label{#1}}

\newcommand{\e}{\operatorname{e}}

\newcommand{\length}{\operatorname{\lambda}}

\newcommand{\Ass}{\operatorname{Ass}}
\newcommand{\Assh}{\operatorname{Assh}}
\newcommand{\Min}{\operatorname{Min}}

\newcommand{\Ann}{\operatorname{Ann}}

\newcommand{\fm}{{\mathfrak m}}

\newcommand{\fa}{{\mathfrak a}}

\newcommand{\fc}{{\mathfrak c}}

\newcommand{\ringR}{\text{$(R,\fm,k)$ }}

\newcounter{hours}\newcounter{minutes}

\newcommand{\excise}[1]{}

\begin{document}
\title[An inequality involving tight closure and parameter ideals]
{An inequality involving tight closure and parameter ideals}
\author[C.~Ciuperc\u{a}]{C\u{a}t\u{a}lin Ciuperc\u{a}}
\author[F.~Enescu]{Florian Enescu}
\address{Department of Mathematics, University of California, Riverside,
CA 92521 USA}
\email{ciuperca@math.ucr.edu}
\address{Department of Mathematics, University of Utah, Salt Lake City,
UT  84112 USA and The Institute of Mathematics of the Romanian Academy, Bucharest, Romania}
\email{enescu@math.utah.edu}

\maketitle

\begin{center}
{\it Dedicated to the memory of Professor Nicolae Radu}
\end{center}

\begin{abstract} We establish an inequality involving colengths of the tight closure of ideals of systems of parameters in local rings with some mild conditions. As an application, we prove and refine a result by Goto and Nakamura, conjectured by Watanabe and Yoshida, which states that the Hilbert-Samuel multiplicity of a parameter ideal is greater than  or equal to the colength of the tight closure of the ideal. 
      
\end{abstract}

\bigskip
\section*{Introduction}
In recent years, the Hilbert-Kunz multiplicities have generated constant interest among researchers in commutative algebra. The development of tight closure theory and its connection to Hilbert-Kunz multiplicities has provided a fresh perspective and led to new discoveries. Among them, a characterization of regular local rings in terms of the Hilbert-Kunz multiplicity of the ring at its maximal ideal, due to K.-i.~Watanabe and Y.~Yoshida, stands out. Under mild conditions, they proved that a local ring $\ringR$ is regular if and only if the Hilbert-Kunz multiplicity at $\fm$ equals $1$ \cite[Theorem 1.5]{WY}. A short and elegant proof of this theorem has also been given by C.~Huneke and Y.~Yao (\cite{HY}). In proving their result, Watanabe and Yoshida were led to a conjecture that ties the Hilbert-Samuel multiplicity of a parameter ideal to the colength of the tight closure of that ideal. Before going further, we would like to give precise definitions to the concepts that appear in our discussion.

Let $(R,\fm)$ be a local ring of positive characteristic $p$. If $I$ is an ideal in $R$, then $I^{[q]}=(i^q: i \in I)$, where $q=p^e$ is a power of the characteristic. Let $R^{\circ} = R \setminus \cup P$, where $P$ runs over the set of all minimal primes of $R$. An element $x$ is said to belong to the {\it tight closure} of the ideal $I$ if $cx^q \in I^{[q]}$ for all sufficiently large $q=p^e$. The tight closure of $I$ is denoted by $I^\ast$. By a ${\it parameter \ ideal}$ we mean here an ideal generated by a full system of parameters in $R$. For an $\fm$-primary ideal $I$, one can consider the Hilbert-Samuel multiplicity and the Hilbert-Kunz multiplicity.

\begin{Definition}
 Let $I$ be an $\fm$-primary ideal in $(R,\fm)$.

1. {\it The Hilbert-Samuel multiplicity of $R$ at $I$} is defined by $\e (I)= \e(I,R) := \displaystyle\lim_{n \to \infty} d! \frac{\length(R/I^n)}{n^d}$. The limit exists and it is positive.

2. {\it The Hilbert-Kunz multiplicity of $R$ at $I$} is defined by $\e _{HK} (I)= \e _{HK}(I,R): = \displaystyle\lim_{q \to \infty}  \frac{\length(R/I^{[q]})}{q^d}$.  Monsky has shown that this limit exists and  is positive.
\end{Definition}

It is known that for parameter ideals $I$, one has $\e(I) = \e_{HK}(I)$. The following sequence of inequalities is also known to hold:
$${\rm max} \{ 1, \frac{1}{d!} \e (I) \} \leq \e_{HK} (I) \leq \e(I)$$
for every $\fm$-primary ideal $I$.

Let $\Assh(R)= \{ P \in \Min(R) : {\rm dim} (R) = {\rm dim}(R/P) \}$. Watanabe and Yoshida have shown that whenever $\Ass(\widehat{R}) = \Assh(\widehat{R})$, the local ring $(R,\fm)$ is regular if and only if $\e_{HK}(\fm) =1$ \cite[Theorem 1.5]{WY}. Under the same assumption that $\Ass (\widehat{R}) =\Assh (\widehat{R})$, they conjectured  that when $I$ is a parameter ideal one has $\e (I) \geq \length(R/I^\ast)$ and that the equality occurs for one (and hence all) parameter ideals if and only if $R$ is Cohen-Macaulay and $F$-rational \cite[Conjecture 1.6]{WY}. Watanabe and Yoshida have also given an affirmative answer to the conjecture for some particular classes of parameter ideals. The conjecture has been proven by S.~Goto and Y.~Nakamura under very mild assumptions \cite[Theorems 1.1 and 1.2]{GN}.

\begin{Theorem}[Goto-Nakamura] 
\label{goto-nak}
Let $(R,\fm)$ be a homomorphic image of a Cohen-Macaulay local ring of characteristic $p>0$.
\begin{enumerate}[1.]
\item Assume that $R$ is equidimensional. Then $\e (I) \geq \length(R/I^\ast)$ for every parameter ideal $I$. 
\item Assume that $\Ass (R) = \Assh (R)$. If $\e (I) = \length(R/I^\ast)$ for some parameter ideal $I$, then $R$ is a Cohen-Macaulay $F$-rational local ring.
\end{enumerate}
\end{Theorem}

In their proof they employ the notion of filter regular sequences. In fact, their proof of part 2 of the Theorem~\ref{goto-nak} is intricate and uses sequences with a property stronger than that of filter regular sequences. Their aim is to reduce the problem to the case of an FLC ring of dimension 2. 

Our main result, Theorem~\ref{prop-1} in Section~1, is an inequality that involves colengths 
of a certain family of parameter ideals. As an immediate application, we refine part 1 of the result of Goto and Nakamura  (as in Remark~\ref{strong-gn}). We also provide a short proof of part 2 of their theorem under some mild conditions.

We would like to thank Craig Huneke for comments that allowed us to improve the manuscript. In particular, Remark~\ref{rem-hun} was suggested by him.

\section{The main result}

 The following theorem is the core of this note.
\begin{Theorem}\label{prop-1}
Let $(R,\fm)$ be an equidimensional local Noetherian ring of characteristic $p > 0$ which is a homomorphic image of a Cohen Macaulay ring and let $x_1, x_2, \ldots, x_d$ be a system of parameters. Then for every $k_1, k_2, \ldots, k_d \geq 1$ we have
\begin{equation}\label{eq-1}
\length\big(R/(x_1^{k_1}, x_2^{k_2}, \ldots, x_d^{k_d})^\ast\big) \geq k_1 k_2\ldots k_d \length\big(R/(x_1, x_2, \ldots, x_d)^\ast\big).
\end{equation}
\end{Theorem}

First, we would like to state the following: 
\begin{Lemma}\label{colon}
Let $(R,\fm)$ be  a local equidimensional ring which is a homomorphic image of a Cohen-Macaulay ring. Let $x_1,\ldots,x_d$ be a system of parameters in $R$. Then 
$$(x_1^{t+s},x_2,\ldots,x_d)^\ast : x_1^{t} \subset (x_1^s,x_2,\ldots,x_d)^\ast$$ for every positive integers t and s.
\end{Lemma}
\begin{proof}
The argument is given essentially in~\cite[(2.3)]{H} and  is also implicit in \cite[(7.9)]{HH}.
\end{proof}

\begin{proof}[Proof of Theorem 1.1]
It is enough to prove that for every system of parameters $x_1, x_2, \ldots, x_d$ and every $k \geq 1$ we  have 
$$\length\big(R/(x_1^{k}, x_2, \ldots, x_d)^\ast\big) \geq k \length\big(R/(x_1, x_2, \ldots, x_d)^\ast\big).$$

Denote $J=(x_1,x_2,\ldots, x_d)$ and set $m=\length(R/J^\ast)$.  Take a filtration of $J^\ast \subseteq R$,
$$J^\ast=\fa_0 \subseteq \fa_1\subseteq  \ldots \subseteq \fa_{m-1} \subseteq \fa_m=R,$$
such that $\length(\fa_{i}/\fa_{i-1})=1$ and write $\fa_i=\fa_{i-1}+(y_i)$ with $y_i \notin a_{i-1}$ and $y_i\fm \subseteq 
\fa_i$.
For each  $n \in \{1,\ldots,k \}$, let $L_i=(x_1^{n}, x_2, \ldots, x_d)^\ast + x_1^{n-1}(y_1,\ldots, y_i)$. Note that $L_m= (x_1^{n}, x_2, \ldots, x_d)^\ast + (x_1^{n-1})$, as $y_m$ is a unit in $R$. 
Then consider  the following filtration:
\begin{equation}\label{eq-2}
(x_1^{n}, x_2, \ldots, x_d)^\ast=L_0 \subseteq L_1 \subseteq \ldots \subseteq L_m \subseteq (x_1^{n-1}, x_2, \ldots, x_d)^\ast. 
\end{equation}

We claim that $\length(L_i/L_{i-1})=1$ for every $i \in \{1,\ldots, m \}$. Indeed, since $y_i\fm \subseteq 
\fa_i=J^\ast +(y_1,\ldots,y_i)$, we have 
$$\fm x_1^{n-1}y_i \subseteq x_1^{n-1}\big(J^\ast +(y_1,\ldots,y_i)\big) \subseteq (x_1^n,x_2,\ldots,x_d)^\ast + x_1^{n-1}(y_1,\ldots,y_i)=L_{i-1}.$$ 
Then  $\fm L_i=\fm(L_{i-1}+ x_1^{n-1}y_i) \subseteq L_{i-1}$, hence $\length(L_i/L_{i-1})\leq 1$. 
On the other hand, we also have $L_i \neq L_{i-1}$. If not, then 
$$x_1^{n-1}y_i \in L_{i-1}=(x_1^n,x_2,\ldots,x_d)^\ast + x_1^{n-1}(y_1,\ldots,y_{i-1}),$$ 
so there exist $r_1,\ldots,r_d$ such that 
$$x_1^{n-1} (y_i - \sum_{j=1}^{i-1}r_jy_j) \in (x_1^n,x_2,\ldots,x_d)^\ast.$$
Since $x_1,\ldots, x_d$ is a system of parameters, by the previous Lemma it follows that 
$$y_i - \sum_{j=1}^{i-1} r_j y_j \in (x_1^n,x_2,\ldots,x_d)^\ast: x_1^{n-1} \subseteq  (x_1,x_2,\ldots,x_d)^\ast=J^\ast,$$  hence $y_i \in J^\ast +(y_1,\ldots,y_{i-1})=\fa_{i-1}$. This contradicts the choice of $y_i$  and the claim is proved.

From (\ref{eq-2}) we then obtain
\begin{equation}\label{eq-3}
\length\big((x_1^{n-1}, x_2, \ldots, x_d)^\ast/(x_1^{n}, x_2, \ldots, x_d)^\ast\big) \geq  m=\length(R/(x_1,x_2,\ldots,x_d)^\ast),
\end{equation} 
with equality if and only if $L_m=(x_1^{n-1}, x_2, \ldots, x_d)^\ast$, that is
\begin{equation}\label{eq-4}
 (x_1^{n-1}, x_2, \ldots, x_d)^\ast=(x_1^{n}, x_2, \ldots, x_d)^\ast +(x_1^{n-1}).
\end{equation}
Finally, as (\ref{eq-3}) holds for every $n \in \{1,\ldots,k\}$,  we get 
$$\length\big(R/(x_1^{k}, x_2, \ldots, x_d)^\ast\big) \geq k \length\big(R/(x_1, x_2, \ldots, x_d)^\ast\big),$$
and the proof is finished.
\end{proof}

\begin{Remark}\label{rem-prop-1} Under the same assumptions, if $$\length\big(R/(x_1^{k}, x_2, \ldots, x_d)^\ast\big) = k \length\big(R/(x_1, x_2, \ldots, x_d)^\ast\big),$$
then 
$$(x_1,x_2,\ldots,x_d)^\ast=(x_1,x_2,\ldots,x_d)+ (x_1^k,x_2,\ldots,x_d)^\ast.$$
Indeed, by (\ref{eq-3}) the first equality implies that 
$$\length\big((x_1^{n-1}, x_2, \ldots, x_d)^\ast/(x_1^{n}, x_2, \ldots, x_d)^\ast\big) = \length(R/(x_1,x_2,\ldots,x_d)^\ast)$$
for every $n \in \{1,\ldots,k \}$. By (\ref{eq-4}) we then obtain
\begin{equation}\label{eq-6}
(x_1^{n-1}, x_2, \ldots, x_d)^\ast=(x_1^{n}, x_2, \ldots, x_d)^\ast +(x_1^{n-1})\quad \text{for } 1 \leq n \leq k,
\end{equation}
which by iteration yields
$$(x_1,x_2,\ldots,x_d)^\ast=(x_1,x_2,\ldots,x_d)+ (x_1^k,x_2,\ldots,x_d)^\ast.$$ 
\end{Remark}
\begin{Corollary}\label{cor-1} Let $(R,\fm)$ be an equidimensional local Noetherian ring of characteristic $p > 0$ which is a homomorphic image of a Cohen Macaulay ring and let $I$ be a parameter ideal. Then  
$$\length\big(R/(I^{[p]})^\ast \big) \geq p^d \length(R/I^\ast).$$
\end{Corollary}

\begin{Remark}\label{rem-1} For every $q=p^e$ set $$a_e= \length\big(R/(I^{[q]})^\ast \big)/q^d \quad\text{and}\quad b_e=\length(R/I^{[q]})/q^d.$$ 
If $(R,\fm)$ is equidimensional and homomorphic image of a Cohen-Macaulay ring and 
$I$ is a parameter ideal, Corollary~\ref{cor-1} shows that $\{a_e\}_{e \geq 0}$ is an increasing sequence. On the other hand, it is  known that  $\{b_e\}_{e \geq 0}$ is a decreasing sequence whose limit is the Hilbert-Kunz multiplicity $\e_{HK}(I)$. In our case  $I$ is a parameter ideal, so $\e_{HK}(I)= \e(I)$.
It is also clear that for every $e \geq 0$ we have
\begin{equation}\label{eq-5}
\length(R/I^\ast) \leq a_e \leq b_e \leq \length(R/I).
\end{equation}
\end{Remark}
\begin{Remark}Assume   that $(R,\fm)$   has a test element  and let $I$ be an $\fm$-primary ideal in $R$. Then $\lim_{q \to \infty} \length\big((I^{[q]})^\ast /I^{[q]}) \big)/q^{d}=0$. Indeed, fix $c$ a test element. Then $(I^{[q]})^\ast \subseteq (I^{[q]} : c)$ so it is enough to show that $\lim_{q \to \infty} \length\big((I^{[q]}:c)/I^{[q]}) \big)/q^{d}=0$. But $\length\big((I^{[q]}:c)/I^{[q]}) \big) = \length \big(R/(c, I^{[q]}) \big)$ and $\lim_{q \to \infty} \length (\big(R/(c, I^{[q]}) \big)/q^d =0$ because $R/cR$ is $(d-1)$-dimensional and $\lim_{q \to \infty} \length (\big(R/(c, I^{[q]}) \big)/q^{d-1}$ represents the Hilbert-Kunz multiplicity of the image of $I$ in $R/cR$.

This also shows that $\lim_{q \to \infty} \length\big(R/(I^{[q]})^\ast \big)/q^d$
exists and equals $e_{HK}(I)$ for an $\fm$-primary ideal $I$. In particular, for a parameter ideal $I$,  
$\lim_{q \to \infty} \length\big(R/(I^{[q]})^\ast \big)/q^d =e(I)$, since 
$\e(I)=\e_{HK}(I)$. So, according to our main result, in this case the sequence $a_e=\length\big(R/(I^{[q]})^\ast \big)/q^d$ is increasing and its limit equals $\e(I)$.
\end{Remark}
As an  immediate consequence of the  Remark~\ref{rem-1}, we obtain  the following result proved by  Goto and Nakamura \cite[Theorem 1.1]{GN}. It had been conjectured and proved in some special cases by Watanabe and Yoshida \cite{WY}.
\begin{Corollary}[Goto--Nakamura]Let $(R,\fm)$ be an equidimensional local Noetherian ring of characteristic $p > 0$ which is a homomorphic image of a Cohen Macaulay ring and let $I$ be a parameter ideal. Then 
$$\e(I) \geq \length(R/I^\ast).$$
\end{Corollary}
\begin{proof} From (\ref{eq-5}) we get $\displaystyle\length(R/I^\ast) \leq \lim_{e \to \infty}b_e= \e_{HK}(I)=\e(I)$.
\end{proof}

\begin{Remark}
\label{strong-gn}
In fact, under the conditions of the above Corollary, one has that

$$\e(I) \geq \length\big(R/(I^{[q]})^\ast \big)/q^d \geq \length(R/I^\ast).$$
\end{Remark}

\begin{proof}
On one hand, $q^d \e(I) = \e(I^{[q]}) \geq \length\big (R/(I^{[q]})^\ast \big)$. On the other hand, we have that $\length\big(R/(I^{[q]})^\ast \big)/q^d \geq \length(R/I^\ast)$, according to the Theorem~\ref{prop-1}.
\end{proof}

As mentioned in the introduction, Watanabe and Yoshida \cite{WY} also conjectured that if $R$ is unmixed $(\Ass (\widehat{R}) = \Assh(\widehat{R}) )$ and $\e(I)=\length(R/I^\ast)$ for some parameter ideal $I$, then $R$ is a Cohen-Macaulay  F-rational ring. If $R$ is a homomorphic image of a Cohen-Macaulay ring and $\Ass (R) =\Assh (R)$, Goto and Nakamura \cite[Theorem 1.2]{GN} proved that the conjecture holds true. Using some considerations on the line of the argument employed in  Theorem~\ref{prop-1}, we are able to give a much shorter proof of this result in the case when either $\widehat{R}$ is reduced or $R$  has a parameter test element.

\begin{Corollary}[Goto--Nakamura]Let $(R,\fm)$ be an equidimensional local Noetherian ring of characteristic $p > 0$ which is a homomorphic image of a Cohen Macaulay ring. Assume that either $\widehat{R}$ is reduced or that $R$ has a (parameter) test element and ${\rm Ass}(R)={\rm Assh}(R)$. If $\e(I)=\length(R/I^\ast)$ for some parameter ideal $I$, then $R$ is a Cohen-Macaulay F-rational ring.
\end{Corollary}
\begin{proof} Let $I=(x_1,\ldots,x_d)$. By (\ref{eq-5}), for every $e \geq 0$ we have
$$\length(R/I^\ast) \leq a_e \leq \lim b_n=\e_{HK}(I)=\e(I),$$
hence
$$\length(R/I^\ast)=\length\big(R/(I^{[q]})^\ast\big)/q^d \quad \text{for every } q=p^e.$$
On the other hand, by Proposition~\ref{prop-1}, for  any $i \in \{1,\ldots,d \}$and all $q=p^e$  we  have
$$\length\big(R/(I^{[q]})^\ast\big)  \geq q^{d-i}\length\big(R/(x_1^q,\ldots,x_i^q,x_{i+1}\ldots,x_d)^\ast\big) \geq   q^d \length\big(R/I^\ast\big)=\length\big(R/(I^{[q]})^\ast\big), $$
 which implies that
$$\length\big(R/(x_1^q,\ldots,x_i^q,x_{i+1},\ldots,x_d)^\ast\big) = q\length\big(R/(x_1^q,\ldots,x_{i-1}^q,x_i,\ldots,x_d)^\ast\big).$$ 

Applying successively  Remark~\ref{rem-prop-1}, we  obtain
\begin{align*}
(x_1,x_2,\ldots,x_d)^\ast &=(x_1,x_2,\ldots,x_d)+ (x_1^q,x_2,\ldots,x_d)^\ast\\
&=(x_1,x_2,\ldots,x_d)+ (x_1^q,x_2^q,x_3\ldots,x_d)^\ast\\
&\ldots\\
&=(x_1,x_2,\ldots,x_d)+ (x_1^q,x_2^q,\ldots,x_d^q)^\ast,
\end{align*}
so $I^\ast=I+(I^{[q]})^\ast$ for all $q=p^e$. 
In particular, $I^\ast \subseteq I + \overline{\fm^q}$ for all $q=p^e$, or equivalently,
$I^\ast \subseteq I + \overline{\fm^n}$ for all $n \geq 1$.

First we consider the case when $\widehat{R}$ is reduced. Passing to the completion $\widehat{R}$, we get  $I^\ast\widehat{R} \subseteq I\widehat{R} + \overline{\fm^n \widehat{R}}$ for all $n \geq 1$. Since $\widehat{R}$ is reduced,   $\bigcap_{n}\overline{\fm^n \widehat{R}}=0$ and by  Chevalley's Lemma, for each $k \geq 1$ there exists $n_k$ with $\overline{\fm^{n_k}\widehat{R}} \subseteq \fm^k\widehat{R}$. This implies that $I^\ast\widehat{R} \subseteq I\widehat{R} + \fm^k\widehat{R}$ for all $k$, hence $I^\ast\widehat{R} \subseteq I\widehat{R}$, or equivalently, $I^\ast=I$. As $I$ is a parameter ideal, this shows that $R$ is a Cohen-Macaulay F-rational ring.

Now, assume that $R$ admits a (parameter) test element $c \in R^{\circ}$. Then $(I^{[q]})^\ast \subseteq (I^{[q]} :c)$. This shows that $c (I^{[q]})^\ast  \subseteq Rc \cap I^{[q]}$. By the Artin-Rees Lemma, one can find $k$ such that $ Rc \cap I^{[q]} \subseteq cI^{q-k}$ for all sufficiently large $q$. In conclusion, $(I^{[q]})^\ast \subseteq (0:c)+ I^{q-k}$. In the first part of proof, we have seen that equality stated in the hypothesis implies that
$I^\ast=I+(I^{[q]})^\ast$ for all $q=p^e$. Hence, $I^\ast \subseteq I + (0:c)+ I^{q-k}$ and taking the intersection over all sufficiently large $q$ we get that $I^\ast \subseteq I +(0:c)$. One can notice now that whenever $\Ass(R) =\Assh(R)$, $c \in R^{\circ}$ implies that $c$ is a nonzerodivisor on $R$, so $(0:c)=0$. In conclusion, $I$ is tightly closed and therefore $R$ is $F$-rational and Cohen-Macaulay.
\end{proof}
\begin{Remark}\label{rem-hun}
The assumption that $R$ admits a parameter test element that appears
in one of the cases of the above Corollary is a mild condition which is
generally easy to test in practice. In fact, if one employs the notion of 
limit closure instead of tight closure, the assumption can
be removed when $R$ is a homomorphic image of a Gorenstein ring and $\Ass(R)= \Assh(R)$.

For any parameter ideal $I$ generated by a system of parameters
$x_1,\ldots,x_d$, one can define the {\it limit closure} of
$I=(x_1,\ldots,x_d)$ by $I^{lim}:= \bigcup_t (x_1^{t+1},\ldots,x_d^{t+1}):(x_1\cdots x_d)^t$.
In general, $I \subseteq I^{lim}
\subseteq I^\ast$ if $R$ is equidimensional and homomorphic image of a Cohen-Macaulay ring. (For more details on the limit closure we refer the reader to~\cite{H,HS}.) The arguments presented in the paper use only the
``colon-capturing'' part of the tight closure and in fact work for the
limit closure too. Hence, one can obtain similar inequalities as in our Theorem that are valid for the limit closure operation. 

Returning to the second
part of Corollary, whenever $R$ is a homomorphic image of a Gorenstein ring
one can prove directly the existence of an element $c$ that multiplies
$(I^{[q]})^{lim}$ into $I^{[q]}$. Indeed, if one denotes $\fa _i = \Ann_R(H^i_{\fm}(R))$, it is known that the ideal $\fc = \fa_0 \cdots \fa_{d-1}$ kills all the modules $(x_1^t,\ldots,x_k^t):x_{k+1}^t/(x_1^t,\ldots,x_k^t)$ for all positive integers $k \leq d-1$ and all positive integers $t$. Then, as in~\cite[p.~208]{S},  one can show that $\fc ^d$ multiplies $(I^{[q]})^{lim}$ into $I^{[q]}$. Since $\dim R/\fa_i \leq i$ for all $i$ \cite[8.1.1(b)]{BH} and  $\Ass(R)= \Assh(R)$,  the ideal $\fc ^d$  contains a nonzerodivisor $c$ which, therefore,  has the above stated property.

Going back to the Corollary, the inequalities $\e(I) \geq \length ( R/I^{lim}) \geq
\length (R/I^\ast)$ coupled with  the hypothesis  $\e(I)=\length(R/I^\ast)$ give that $I^{lim} =I^{\ast}$. As in the proof of the Corollary, one can now conclude that $I^{lim}=I+(I^{[q]})^{lim}$ for all $q=p^e$. The existence of a nonzerodivisor $c$ as above allows us to continue the proof in similar fashion and conclude that $I = I^{lim}$.  So, $I = I^{lim} =I^\ast$ and therefore $R$ is $F$-rational and Cohen-Macaulay.

The reader should also note that one good property of the limit closure is that it commutes
with the completion at the maximal ideal. The same is true for the multiplicity of a parameter ideal. Consequently, if one chooses to pass to the completion first, the completion is a homomorphic image of a regular ring and  the proof outlined above will work similarly if one assumes now that $R$ is a homomorphic image of a Cohen-Macaulay ring and $\Ass(\widehat{R})= \Assh(\widehat{R})$. 
\end{Remark}

The previous results can be easily extended to the class of arbitrary
                 ideals primary to the maximal ideal.

\begin{Corollary}
Let $\ringR$ be an equidimensional local Noetherian ring of characteristic $p > 0$ which is a homomorphic image of a Cohen Macaulay ring. Assume that either $R$ and $\widehat{R}$ have a common test element, or $k$ is infinite.
\begin{enumerate}[1.]
\item For every $\fm$-primary ideal $I$, we have $$\e(I) \geq \length(R/I^\ast).$$
\item If $\e(I)=\length(R/I^\ast)$ for some $\fm$-primary ideal $I$, then $R$ is a Cohen-Macaulay F-rational ring.
\end{enumerate}
\end{Corollary}
\begin{proof}
If $k$ is finite and $R$ and $\widehat{R}$ have a common test element, we can enlarge the residue field of $\ringR$ such that $k$ is an infinite field, by passing to $R[X]_{\fm R[x]}$. (The tight closure of $I$ commutes with this base change by~\cite[Theorem 7.16]{HH2}.) Hence we can assume that the residue field of $R$ is infinite. This allows us to consider a reduction $J$ of $I$ such that $J$ is a parameter ideal. Then $\e (I)= \e(J) \geq \length(R/J^\ast) \geq \length(R/I^\ast)$, where the last inequality holds because $J \subseteq I$ and hence $J^\ast \subseteq I^\ast$.

To prove part 2, notice that the equality for $I$ implies that  $\e(J) \geq \length(R/J^\ast)$, for $J$ chosen as above. Hence $R$ is Cohen-Macaulay and $F$-rational by Theorem~\ref{goto-nak}.
\end{proof}

\end{document}